\documentclass[a4paper]{amsart}

\usepackage{amssymb,amsmath,latexsym,amsthm,amsrefs,graphicx,array,float,stmaryrd,empheq,enumitem,url}

\usepackage{amsthm}
 
 \newtheorem{theorem}{Theorem}
 \newtheorem{remark}{Remark}
 
 \newtheorem{lemma}{Lemma}
 \newtheorem{conjecture}{Conjecture}
 
\theoremstyle{definition}
\newtheorem{definition}{Definition}[section]

% ...

\begin{document}

\title
[Remarks on an interpolation between Wilson's theorem and Giuga's conjecture]
{Remarks on an interpolation between Wilson's theorem and Giuga's conjecture}

\author{Thomas Sauvaget}


\email{thomasfsauvaget@gmail.com}


\thanks{This version of the manuscript takes into account helpful comments and encouragements provided by several  referees of previous versions, to whom the author is very grateful. The author also wishes to thank all those who contributed, directly or indirectly, to make Sage and its libraries such a useful free software.}

  



\begin{abstract}
A family of congruences interpolating between those of Wilson and Giuga is constructed. Several elementary results are presented, in order to present a possible approach to establishing Giuga's conjecture. 
\end{abstract}

\maketitle

\section{Introduction}

Let $p$ be a prime, and $j$ a positive integer such that $1\leq j\leq p-1$. 
Then by Fermat's Little Theorem we have $j^{p-1}\equiv 1 \pmod{p}$. 
And thus $\sum_{j=1}^{p-1}j^{p-1}\equiv p-1 \equiv -1 \pmod{p}$, 
an observation of Sierpi\'nski \cite{G}. 
\bigskip

In 1950, Giuga \cite{GG} asked whether the converse holds, 
i.e. if any integer $n$ such that 
$\sum_{j=1}^{n-1}j^{n-1}\equiv -1 \pmod{n}$ is necessarily prime. 
He provided support for that conjecture by showing that any counterexample, 
i.e. any composite integer satisfying the congruence, 
is greater than $10^{1000}$. 
This lower bound was later improved many times over the ensuing decades, by numerical means : to $10^{1700}$ by Bedocchi \cite{B}, to $10^{13887}$ by Borwein, Borwein, Borwein \& Girgensohn \cite{BBBG}, and finally to $10^{19908}$ by Borwein, Maitland \& Skerritt \cite{BMS}. These works nevertheless have 
arithmetical limitations, in that the method used intrinsically cannot treat any number having more than 8,135 prime factors.
\bigskip

On the other hand, Tipu \cite{T} proved by analytic means that for any real number $x$, the number of counterexamples to Giuga's conjecture 
$G(x):=\#\{ n<x: n \text{ is composite and }  
\sum_{a=1}^{n-1}a^{n-1}\equiv -1 \pmod{n}\}$, if any,  
is at best of moderate growth: $G(x)\ll \sqrt{x}\log x$. This was later 
improved by Luca, Pomerance and Shparlinski \cite{LPS} 
to $G(x)=o(\sqrt{x})$.  
\bigskip

Agoh \cite{A} showed in 1995 that conjecturing that $n$ is prime iff $nB_{n-1}\equiv -1 \pmod{n}$ (where the $(B_{n})_{n\geq 0}$ are Bernoulli numbers, see e.g. Granville \cite{GR} and Sun \cite{ZWS})  is equi\-valent to Giuga's conjecture (see also Kellner \cite {K}). 
\bigskip

The paper is organized as follows. In the next section, we introduce a family of congruences that interpolates between Wilson's and Giuga's, and prove an important lemma. 
Then, in section 3, we show that the family of congruences can be rewriten using  unsigned Stirling numbers of the first kind. Properties of those numbers observed numerically, combined with the important lemma, would then imply that odd squarefree composite integers cannot satisfy Giuga's congruence. 
Also, in an appendix, we present a combinatorial relation which might be of independent interest (which arised in an earlier incorrect attempt at proving Giuga's conjecture).
\bigskip

{\bf Notations.} Let $\mathbb{N}_2:=\mathbb{N}-\{0,1\}$. The set of primes is denoted by $\mathcal{P}:=\{2, 3, 5, 7, 11,\dots\}$. The double-brackets denote 
sets of consecutive integers: $\llbracket a;b \rrbracket := \{ i\in\mathbb{N} | a\leq i \leq b \}$.
Define also for each $k\in\mathbb{N}$ the functions $S_k:\mathbb{N}\rightarrow\mathbb{N}$ by $S_k(n):=\sum_{j=1}^{n-1}j^k$ .

\section{An interpolation of Wilson's and Giuga's congruences}

\subsection{Basic definitions}

\begin{definition}
\label{wilsonfunction}
Let  $ f_W :  \mathbb{N}_2\mapsto \mathbb{N} \text{ be defined by }  f_W(n):=1+(n-1)!$
\end{definition}

Then Wilson's theorem can be recast as:

\begin{theorem} $\forall n\in \mathbb{N}_2 \text{ we have: } 
n\in\mathcal{P} \Leftrightarrow n\mid f_W(n)$.
\end{theorem}

\begin{remark}
Wilson's theorem allows to change the definition of primality, which is negative in nature ($n$ is prime iff it is not the product of strictly smaller integers), into a positive one ($n$ is prime iff it divides some specific larger integer). 
This crucial reversal is at the heart of the philosophy of this paper.
\end{remark}

\begin{definition}
\label{giugafunction}
Let  $f_G :  \mathbb{N}_2\mapsto \mathbb{N} \text{ be defined by } f_G(n):=1+S_{n-1}(n)$.
\end{definition}

Then Giuga's conjecture is:

\begin{conjecture}
\label{giugasconjecture}
$\forall n\in \mathbb{N}_2 \text{ we have: } 
n\in\mathcal{P} \Leftrightarrow n\mid f_G(n)$.
\end{conjecture}

\subsection{An interpolation family}

We introduce the following integer-valued functions. 

\begin{definition}
Let 
$
H_k : \mathbb{N}_2\mapsto \mathbb{N} \text{ be defined by }
H_k(n):=
1+\sum_{i=1}^{n-1}\frac{(n+i-1-k)!}{i!}i^k 
\text{ for } k \in\llbracket 1;n-1 \rrbracket
$
\end{definition}

\begin{remark}
For any given $n\in\mathbb{N}_2$, we have that for any $k\in\llbracket 1;n-1 \rrbracket$ the function $H_k$ is integer-valued, since for each $i\in \llbracket 1;n-1 \rrbracket$ we have $i\leq n+i-1-k$.
\end{remark}

\begin{remark}
We directly obtain $H_{n-1}=f_G$. The relationship between $H_1$ and $f_W$ is described in the proof of the next lemma.
\end{remark}

\begin{lemma}
\label{proofH1}
$\forall n\in \mathbb{N}_2 \text{ we have: } 
n\in\mathcal{P} \Leftrightarrow n\mid H_1(n)$.
\end{lemma}
\begin{proof}
Since by definition $H_1(n)=1+(n-1)!+\frac{n!}{2!}2+\frac{(n+1)!}{3!}3+\cdots +\frac{(2n-3)!}{(n-1)!}(n-1)=1+(n-1)!+n[ \frac{((n-1)!}{2!}2+\frac{(n+1)!}{3!n}3+\cdots +\frac{(2n-3)!}{(n)!}(n-1)]\equiv 1+(n-1)!\equiv f_W(n)\pmod{n}$. The term in the bracket is necessarily an integer since for any $j\in\llbracket1,n-3\rrbracket$ we have that 
$(n+j)!=\underbrace{1\times 2\times \cdots\times (j+2)}_{=(j+2)!} \times (j+3)\times\cdots\times n\times (n+1)\times\cdots\times (n+j) $, so this is simultaneously divisible by $(j+2)!$ and by $n$. 
The result then follows from Wilson's theorem.
\end{proof}
\bigskip

Thus $H_1$ also characterizes the primes. To motivate the strategy employed in the next sections, consider the residues of the values of $H_k(n)$ modulo $n$,  for $k\in\llbracket 1;34 \rrbracket$  and $n\in\llbracket 5;35 \rrbracket$, which are provided in the following table:

\begin{figure}[H]
\label{numdata}
\centering
\includegraphics[scale=0.44]{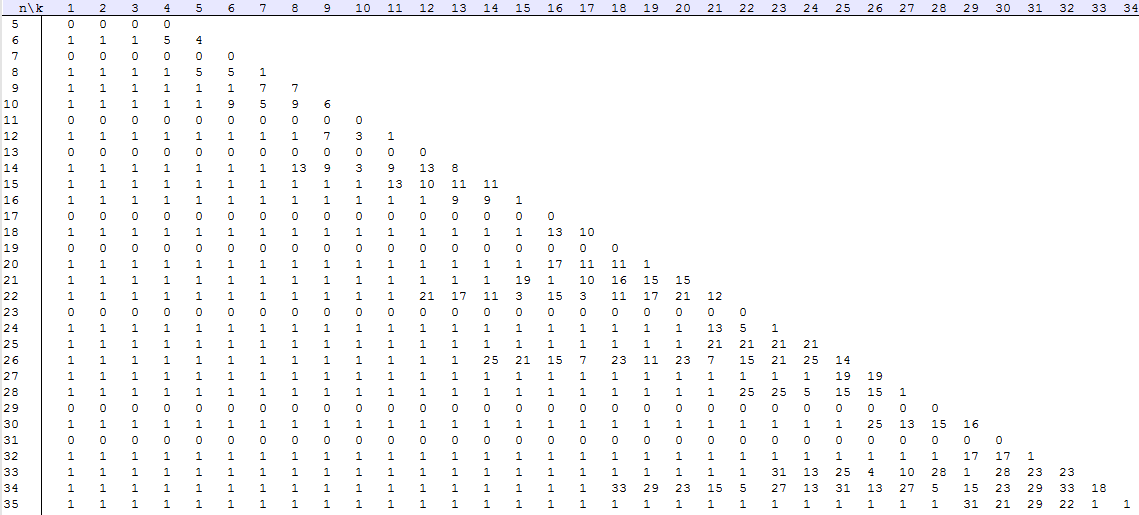}
\caption{Residues of the values of $H_k(n)$ modulo $n$.}
\end{figure}

On this very limited data set, the following observations can be made: 
\begin{itemize}
\item for any $k\in \llbracket 1;n-1 \rrbracket$, $H_k(n)\equiv 0\pmod{n}$ for all $n$ prime
\item for any $k\in \llbracket 1;n-1 \rrbracket$, $H_k(n)\not\equiv 0\pmod{n}$ for any composite $n$
\item in the case $n$ composite, $H_k(n)\equiv 1\pmod{n}$ for $k\leq \frac{n}{2}$
\end{itemize}
\bigskip

The last point can be easily established beyond $n=35$ for squarefree composites by the following lemma:

\begin{lemma}
\label{importantLemma}
Let $n\in\mathbb{N}_2$ be composite and squarefree. Then  $H_k(n)\equiv 1\pmod{n}$ for $k\leq \frac{n}{2}$.
\end{lemma}

\begin{proof}
Given $n$ composite, let its prime factorization be 
$n=\prod_{j=1}^J q_j$, 
where $q_1<\dots <q_J$ are distinct primes. We necessarily have that 
$2\leq q_1<\dots <q_J \leq \frac{n}{2}$.
\medskip

Let $k\in \llbracket 1;\lfloor\frac{n}{2}\rfloor \rrbracket$. So, for each $i\in \llbracket 1;n-1 \rrbracket$ we have  $n+i-1-k\geq \lfloor\frac{n}{2}\rfloor +i-1\geq q_J +i-1\geq q_J>\cdots >q_1\ \ (\diamond )$.
\medskip

Recall the definition of $H_k(n)$. Since $\frac{(n+i-1-k)!}{i!}$ is a product of $n+i-1-k$ consecutive integers starting at $i+1$, each term $\frac{(n+i-1-k)!}{i!}i^k$ contains a product of $n+i-k$ consecutive integers starting at $i$.
\medskip

So, given the inequalities $ (\diamond )$, we infer that for each $i\in\llbracket 1;n-1 \rrbracket$, each of the $q_j$ 
is a factor of at least one of those consecutive integers,  hence each term 
$\frac{(n+i-1-k)!}{i!}i^k$ is a multiple of $n$. So only the first term $1$ in $H_k(n)$ remains modulo $n$.
\end{proof}
\bigskip

Going beyond $k=\lfloor\frac{n}{2}\rfloor$ would require a more detailed analysis of the structure of those values of $H_k(n)$. But, given the first and second observations, 
it appears more straighforward, to make the following conjecture: it captures the primality condition (no composite is a false positive), while not entering into the arithmetical details of the values of the residues for $n$ composite.

\begin{conjecture}
\label{conjectureH}
For each $k\in  \mathbb{N}_2 \text{ we have: } 
[n\in\mathcal{P},\text{ with } n\geq k+1] \Leftrightarrow n\mid H_k(n)$.
\end{conjecture}

\begin{remark}
The case $k=n-1$ is Giuga's conjecture. 
\end{remark}

In the next section, we show how lemma \ref{importantLemma} could make it impossible for odd squarefree composite integers to satisfy Giuga's congruence.

\section{Towards a path for the proof of Giuga's conjecture}

\subsection{The case of $n=p$ prime of conjecture \ref{conjectureH}}

We start with the following elementary lemma, whose proof is taken almost verbatim from Nicolas \cite{NMSE}.

\begin{lemma}
\label{lemmaSk}
Let $p$ be an odd prime. Then for any $k\in\mathbb{N}^*$ such that $\gcd(p-1,k)=1$ we have $S_k(p)\equiv 0\pmod{p}$.
\end{lemma}

\begin{proof}
Let $b$ be any number relatively prime to $p$. Recall that the numbers
 $b, 2b, 3b, $ 
 
 $4b,\dots , (p-1)b$ are a reduced residue class modulo $p$. 
\medskip

It follows that
$b^k+(2b)^k+(3b)^k+\dots+((p-1)b)^k\equiv 1^k+2^k+3^k+\dots+(p-1)^k\pmod{p} \ \ (*)$.
But $(ib)^k=b^ki^k$, so we can rewrite $(*)$ as
$(b^k-1)(1^k+2^k+3^k+\dots+(p-1)^k)\equiv 0\pmod{p}\ \ (**)$.

Let $b$ have order $p-1$ modulo $p$. So we are letting $b$ be a primitive root of $p$. Since $p-1$ does not divide $k$, we have $b^k\not\equiv1\pmod{p}$. Then it follows from $(*)$ that $1^k+2^k+3^k+\dots+(p-1)^k\equiv 0\pmod{p}$ .
\end{proof}
\bigskip

Recall the following classical definitions.
\begin{definition}
The unsigned Stirling numbers of the first kind, $\begin{bmatrix}M\\j\end{bmatrix}$, are the coefficients of the expansion of the rising factorial into powers: $\prod_{j=0}^{j=M-1}(x+j) =\sum_{j=0}^{j=M}\begin{bmatrix}M\\j\end{bmatrix}x^j$. 
The Stirling numbers of the first kind, $s(M,j)$,  are the coefficients of the expansion of the falling factorial into powers: $\prod_{j=0}^{j=M-1}(x-j) =\sum_{j=0}^{j=M}s(M,j)x^j$. We have $|s(M,j)|=\begin{bmatrix}M\\j\end{bmatrix}$.
\end{definition}

We can use this and lemma \ref{lemmaSk} to prove the following result, which is one direction of conjecture \ref{conjectureH}.

\begin{theorem}
\label{expansion}
For all odd prime $p$, and all $k\in\llbracket 1;p-2 \rrbracket$  we have $H_k(p)\equiv 0\pmod{p}$.
\end{theorem}
\begin{proof}
Indeed, for any $n$, prime or composite, we can expand $H_k(n)$ as function of $H_{n-1}(n)$, and $\{S_k(n),\dots , S_{n-2}(n)\}$, and $\{\begin{bmatrix}n-k\\1\end{bmatrix},\dots , \begin{bmatrix}n-k\\n-k-1\end{bmatrix}\}$:

\[
\begin{array}{ccl}
H_k(n) &= & 1+\sum_{i=1}^{n-1}\frac{(n+i-1-k)!}{i!}i^k \\
& = & 1+ \sum_{i=1}^{n-1} (n+i-1-k)(n+i-2-k)\dots (1+i) i^k\\
& = & 1+ \sum_{i=1}^{n-1} \left(\sum_{j=0}^{n-k} \begin{bmatrix}n-k\\j\end{bmatrix} i^j\right)i^{k-1}\\
\\
& = & 1+ \sum_{i=1}^{n-1} \left(\begin{bmatrix}n-k\\0\end{bmatrix} i^{k-1}+\dots +\begin{bmatrix}n-k\\n-k\end{bmatrix} i^{n-1}\right)
\end{array}
\]

This can the be rearranged as:
\[
\begin{array}{ccl}
H_k(n) & = & H_{n-1}(n) +  \begin{bmatrix}n-k\\n-k-1\end{bmatrix}S_{n-2}(n)+\dots +\begin{bmatrix}n-k\\1\end{bmatrix} S_{k}(n)
\end{array}
\]
since $\begin{bmatrix}n-k\\0\end{bmatrix}=0$ and $\begin{bmatrix}n-k\\n-k\end{bmatrix}=1$.
\medskip

Now, for $n=p$ prime we mentioned at the begining of this paper that $H_{p-1}(p)\equiv 0\pmod{p}$. Combining this with lemma \ref{lemmaSk} finishes the proof.
\end{proof}
\medskip

We shall not attempt to prove conjecture \ref{conjectureH} in this work. Let us only mention the following simpler conjecture based on numerical data.

\begin{conjecture}
For any odd squarefree composite $n=q_1\cdots q_J$, one has $H_{n-q_J+1}(n)\equiv n+1-q_1\pmod{n}$.
\end{conjecture}

\begin{remark}
This would be, for odd squarefree composites, the first ``non-trivial value of $k$'' i.e. for these $n$ one could show like in lemma 2 that $H_k(n)\equiv 1\pmod{n}$ for $k<n-q_J+1$, while for $k=n-q_J+1$ one would have $H_k(n)$ different from either 0 or 1 modulo $n$.
\end{remark}

\subsection{Steps towards Giuga's conjecture}

Recall from \cite{BBBG} that even, and non-squarefree, composite integers cannot satisfy Giuga's congruence. So in this section we consider only the case of odd squarefree integers, i.e. $n=q_1\dots q_J$ with $3\leq q_1<\dots <q_j\leq \frac{n}{3}$.
\bigskip

The following elementary lemma allows to reduce the problem to the prime factors of $n$. 

\begin{lemma}
Let $n=q_1\cdots q_J$ be odd composite squarefree. Then for any number $M$ we have the equivalence:

\[ 
M\equiv 0\pmod{q_1\cdots q_J} \Leftrightarrow \left\{ 
\begin{array}{l}
M\equiv 0\pmod{q_1}\\  
\vdots \\
M\equiv 0\pmod{q_J}
\end{array}
\right .
\]
\end{lemma}

\begin{proof}
This follows immediately by the Chinese Remainder Theorem.
\end{proof}

We can apply this lemma in two special cases. 
Firstly, the condition $H_k(n)\equiv 0\pmod{n}$ is equivalent to having 

\[
\left\{ 
\begin{array}{l}
H_k(n)\equiv 0\pmod{q_1}\\  
\vdots \\
H_k(n)\equiv 0\pmod{q_J}
\end{array}
\right .
\]
\medskip

So to show $H_k(n)\not\equiv 0\pmod{n}$ it would suffice to show that for at least one $q_j$ we have $H_k(n)\not\equiv 0\pmod{q_j}$.
\bigskip

Secondly, the condition $H_k(n)\equiv 1\pmod{n}$ is equivalent to having 

\[
\left\{ 
\begin{array}{l}
H_k(n)\equiv 1\pmod{q_1}\\  
\vdots \\
H_k(n)\equiv 1\pmod{q_J}
\end{array}
\right .
\]
\medskip

So to show $H_k(n)\not\equiv 1\pmod{n}$ it would suffice to show that for at least one $q_j$ we have $H_k(n)\not\equiv 1\pmod{q_j}$.
\bigskip

Recall the following result on congruences for unsigned Stirling numbers of the first kind, which we take verbatim from Scott \cite{S}.

\begin{lemma}
\label{congruence_stirling}
$p$ is prime $\Leftrightarrow$ for all $k\in\llbracket 2;p-1 \rrbracket$ we have 
$\begin{bmatrix}p\\k\end{bmatrix}\equiv 0\pmod{p}$.
\end{lemma}
\begin{proof}
Work in the polynomial ring $(\mathbb{Z}/p\mathbb{Z})[x]$: 
in that ring $\prod_{i=0}^{p-1}(x-i)=x^p-x$, since both polynomials are monic of degree $p$ and have roots $0,\dots ,p-1$.
\smallskip

Now just equate the coefficients to get the positive part of the result. For the negative result, suppose that $n$ is not prime, let $p$ be a prime dividing $n$, and let 
$m=\frac{n}{p}$. Then

\[
\prod_{i=0}^{n-1}(x-i)=
\prod_{k=0}^{m-1}\prod_{\ell=0}^{p-1}(x-kp-\ell)
=\prod_{k=0}^{m-1}\prod_{\ell=0}^{p-1}(x-\ell)
=\left(\prod_{\ell=0}^{p-1}(x-\ell) \right)^m = (x^p-x)^m=x^m(x^{p-1}-1)^m
\]

\noindent and it follows that $s(n,m)\equiv(-1)^m\not\equiv 0\pmod{p}$. 
Since $\begin{bmatrix}n\\m\end{bmatrix}=|s(n,m)|$ the claim follows.
\end{proof}

The general expansion of theorem \ref{expansion}:
\[
\begin{array}{ccl}
H_k(n) & = & H_{n-1}(n) +  \begin{bmatrix}n-k\\n-k-1\end{bmatrix}S_{n-2}(n)+\dots +\begin{bmatrix}n-k\\1\end{bmatrix} S_{k}(n)
\end{array}
\]

\noindent becomes for $k=q_j$ :
\[
\begin{array}{ccl}
H_{q_j}(n) & = & H_{n-1}(n) +  \begin{bmatrix}n-q_j\\n-q_j-1\end{bmatrix}S_{n-2}(n)+\dots +\begin{bmatrix}n-q_j\\1\end{bmatrix} S_{q_j}(n)
\end{array}
\]

So if we assume that $H_{n-1}(n)\equiv 0\pmod{n}$, which implies  $H_{n-1}(n)\equiv 0\pmod{q_j}$, and combine this with our observation of section 1 that $H_k(n)\equiv 1 \pmod{n}$ for $k\leq \frac{n}{3}$ (since $n$ is odd), so in particular  $H_{q_j}(n)\equiv 1 \pmod{q_j}$, we get:

\[
1\equiv \begin{bmatrix}n-q_j\\n-q_j-1\end{bmatrix}S_{n-2}(n)+\dots +\begin{bmatrix}n-q_j\\1\end{bmatrix} S_{q_j}(n) \pmod{q_j}
\]

\begin{remark}
Numerically, we observe that for the largest prime factor $q_K$ we \emph{always}  seem to have $\begin{bmatrix}n-q_K\\n-q_K-1\end{bmatrix}S_{n-2}(n)+\dots +\begin{bmatrix}n-q_K\\1\end{bmatrix} S_{q_K}(n) \equiv 0 \pmod{q_K}$. If one could prove this, a contradiction would follow, hence one would conclude that one cannot have $H_{n-1}(n)\equiv 0\pmod{n}$, i.e. Giuga's conjecture is true. But the precise way in which that cancellation occurs modulo $q_K$ seems non-trivial, in particular most of the time several terms are non-zero. Note also that we cannot use  lemma \ref{congruence_stirling}, at least not directly, since we are dealing with the composite $n-q_K$.
\end{remark}

\section{Appendix}

In this appendix we present an earlier incorrect attempt of the author, stressing where the error lies. The author is very grateful for the input from multiple referees to finally help him understand the problem.
\bigskip

\subsection{General properties}

\begin{lemma}
\label{ordering}
For each $(n,k)\in\mathbb{N}_2^2$ with $n\geq k+1$, we have :   
$H_{k}(n)<H_{k-1}(n)$.
\end{lemma}
\begin{proof}
With those $n$ and $k$, since $\forall i\in\llbracket 1;n-1 \rrbracket$ 
we have $\frac{(n+i-1-k)!}{i!}i^k = \frac{i}{n+i-1-k} \times \frac{(n+i-k)!}{i!}i^{k-1}$, 
we thus explicitely have 
$H_k(n)=1+\sum_{i=1}^{n-1}\frac{(n+i-1-k)!}{i!}i^k
<1+\sum_{i=1}^{n-1}\frac{(n+i-k)!}{i!} i^{k-1}=H_{k-1}(n)$
\end{proof}

\begin{lemma}
\label{corepart}
For each $(n,k)\in\mathbb{N}_2^2$ with $n\geq k+1$ we have :  
\[
H_k(n)\equiv 1+\sum_{i=1}^k \frac{(n+i-1-k)!}{i!}i^k \pmod{n}
\]
\end{lemma}
\begin{proof}
Immediate.
\end{proof}

\subsection{Passage from $H_k$ to $H_{k+1}$} 

The following result can be established. 

\begin{lemma}
\label{step}
For each $n\in\mathbb{N}_2$ and any $1\leq k\leq n-2$ we have :  
\[
H_{k+1}(n) = H_k(n) - (n-k-1) \sum_{i=1}^{n-1} \frac{(n+i-2-k)!}{i!}i^k
\]
\end{lemma}
\begin{proof}
We have $H_{k+1}(n) =1+ \sum_{i=1}^{n-1} \frac{(n+i-(k+1)-1)!}{i!}i^{k+1} $

$= 
1+ \sum_{i=1}^{n-1} \frac{i}{n+i-k-1}\frac{(n+i-k-1)!}{i!}i^k$, so 
adding and substracting 

$ \sum_{i=1}^{n-1} \frac{n-k-1}{n+i-k-1}\frac{(n+i-k-1)!}{i!}i^k$ to that, one obtains the result.
\end{proof}

The following functions will be useful. 

\begin{definition}
 Let 
 $U_{k} : \mathbb{N}_2\mapsto \mathbb{N} \text{ be defined by }$
 
 \[
U_k(n):=
\sum_{i=1}^{k+1} \frac{(n+i-2-k)!}{i!}i^k \text{ if $k \in\llbracket 1;n-2 \rrbracket$ }
\]
\end{definition}

\begin{remark} For any given $n\in \mathbb{N}_2$ we have that for any 
$k \in\llbracket 1;n-2 \rrbracket$ the quantity $U_{k}(n)$ is an integer, since 
for any $i \in\llbracket 1;k+1 \rrbracket$ we have $i\leq n+i-2-k$.
\end{remark}

\begin{remark} Since $\frac{(n+i-2-k)!}{i!}i^k\equiv 0 \pmod{n}$ for $i\in\llbracket k+2;n-1 \rrbracket$,  we have $H_{k+1}(n)-H_k(n)\equiv (k+1) U_{k}(n)\pmod{n}$.
\end{remark}

\begin{lemma}
For any $n\in\mathbb{N}_2$ we have $U_1(n)\equiv 0 \pmod{n}$. In particular, $H_2$ also characterizes $\mathcal{P}$.
\end{lemma}
\begin{proof}
We have to show that for $n\geq 2$ 
we have $ \frac{(n-2)!}{1!}1^1 + \frac{(n-1)!}{2!}2^1 
\equiv 0  \pmod{n}$. This is true since it simplifies 
to  showing $(n-2)! +(n-1)! \equiv 0  \pmod{n}$, i.e. showing that $(n-2)! [1+(n-1)] \equiv 0  \pmod{n}$, which is true. So for $n\geq 2$ we have $H_2(n)\equiv H_1(n)\pmod{n}$, so lemma \ref{proofH1} implies that $H_2$ also characterizes the  primes. 
\end{proof}
\medskip

\subsection{A combinatorial relationship}

Let us now introduce the following sums.

\begin{definition}
For each $ k\in\mathbb{N}^*$ let 
\[ V_{k}:=\sum_{i=1}^{k} \left [ (-1)^{i-1}  \binom{k+1}{i}   i^k + (-1)^k \binom{k}{i} k^i \right ] 
\]
\end{definition}

\begin{theorem}
\label{keyTheorem} 
For each $ k\in\mathbb{N}^*$ we have $V_{k}=(-1)^{k+1}$.
\end{theorem}

This will be proved as a result of the following two lemmas. The author is very grateful to the referee of a previous version of this work for providing this proof. 
This result might be of independent interest.

\begin{lemma}
\label{lemmaA}
For each $k\in\mathbb{N}^*$ we have $\sum_{i=0}^k {k\choose i} k^i = (k+1)^k$.
\end{lemma}
\begin{proof}
This follows immediately from the binomial theorem.
\end{proof}

\begin{lemma}
\label{lemmaB}
For each $k\in\mathbb{N}^*$ we have 
$\sum_{i=0}^{k+1} (-1)^i {{k+1}\choose i} i^k = 0$.
\end{lemma}
\begin{proof}
Defining for any function $f:\mathbb{Z}\rightarrow\mathbb{Z}$ the forward difference operator $\Delta$ by $\Delta f (j):=f(j+1)-f(j)$, we obtain by induction on $n$ that the $n-$fold iteration of $\Delta$ on $f$ is $\Delta^n f(j)=\sum_{i=0}^n (-1)^{i} {n\choose i} f(j+n-i)$.
\medskip

For $n=k+1$ and the polynomial $f(j)=j^k$, since $f$ is of degree $k<n$ we thus have that $\Delta^n f(j)$ is identically zero. That is, we have 
 $\sum_{i=0}^{k+1} (-1)^i {{k+1}\choose i} (j+k+1-i)^k = 0$. 
\medskip

In particular, when evaluated at $j=-(k+1)$ this 
gives $\sum_{i=0}^{k+1} (-1)^i {{k+1}\choose i} (-i)^k =
(-1)^k \sum_{i=0}^{k+1} (-1)^i {{k+1}\choose i} i^k = 0$, hence the result.
\end{proof}

\bigskip

We now prove theorem \ref{keyTheorem}. Let 
\[
A:=\sum_{i=1}^k {k\choose i} k^i
\]
and 
\[ 
B:=\sum_{i=1}^{k} (-1)^i {{k+1}\choose i} i^k
\]

Then, given the slightly different summation ranges than in the previous lemmas, we get $A= (k+1)^k-1$ and $B=(-1)^k(k+1)^k$. 
So since $(-1)^kA+ B =(-1)^{k+1}$ we obtain $V_k=(-1)^{k+1}$.

\subsection{The error in the current attempt}

Our aim was then to prove the following:

\begin{conjecture}
\label{theoremSum}
With the previous notations, we have for any $n\in\mathbb{N}_2$
and any $1\leq k\leq n-2$ that $U_{k}(n)\equiv 0\pmod{n}$.
\end{conjecture}

Unfortunately it is not true, as was explained to us by (we now realise) several referees. Our aim was to proceed as follows.
\bigskip

The $i=1$ term 
of $U_{k}(n)$ is $(n-1-k)!$. Let us say that $(n-1-k)!\equiv a \pmod{n}$ for some $a\in \llbracket 0;n-1 \rrbracket$. Then for any $j\in \llbracket 1; k \rrbracket$ we have by induction  $(n-1-k+j)!\equiv a \times (n-k)\times \cdots \times (n-1-k+j) \equiv (-1)^j a \frac{k!}{(k-j)!} \pmod{n}$. So $a$ is in factor in each term  of $U_{k}(n)$. 
\medskip

If $a=0$ then $U_{k}(n)\equiv 0\pmod{n}$. If $a\neq 0$, then $U_{k}(n)$ can be factored by $a$, but that is as far as one can proceed. In particular, one \textbf{cannot} write (using the previous paragraph with $i=j+1$) that we have  $U_{k}(n)=a\times\sum_{i=1}^{k+1} \frac{ (-1)^{i-1}\frac{k!}{(k+1-i)!}}{i!}i^{k} = a\times\left(\sum_{i=1}^{k} \frac{1}{k+1}\binom{k+1}{i} (-1)^{i-1} i^{k} + (-1)^{k}\frac{(k+1)^k}{k+1}\right)$, due to the division by $i!$ (non-commutativity of two otherwise true congruences -- this is nothing else than the fact that $H_k(n)\not\equiv 0 \text{ nor } 1\pmod{n}$ as seen in the numerical data). 
\medskip

Had that been the case, recalling that the term in parenthesis is an integer, and multiplying by the quantity $k+1$ (which is  non-zero modulo $n$ since $1\leq k \leq n-2$) we would have obtained  
$(k+1)U_k(n)=a\times\left(\sum_{i=1}^{k} \binom{k+1}{i} (-1)^{i-1} i^{k} + (-1)^{k}(k+1)^k\right)$. Expanding $(k+1)^k$ with the binomial theorem and simplifying, we notice that the resulting expression is $0$ by theorem \ref{keyTheorem}.  So we would have necessarily had $U_k(n)\equiv 0\pmod{n}$...

\bibliographystyle{amsrefs}
\bibliography{biblio_giuga_v3}

\begin{bibdiv}
\begin{biblist}

\bib{A}{article}{
      author={Agoh, T.},
       title={On {G}iuga's conjecture},
        date={1995},
     journal={Manuscripta Math.},
      volume={87},
       pages={501\ndash 510},
}

\bib{B}{article}{
      author={Bedocchi, E.},
       title={Nota ad una congettura sui numeri primi},
        date={1985},
     journal={Riv. Mat. Univ. Parma (4)},
      volume={11},
       pages={229\ndash 236},
}

\bib{BBBG}{article}{
      author={Borwein, D.},
      author={Borwein, J.M.},
      author={Borwein, P.B.},
      author={Girgensohn, R.},
       title={{G}iuga's conjecture on primality},
        date={1996},
     journal={The American Mathematical Monthly},
      volume={103},
       pages={40\ndash 50},
}

\bib{BMS}{article}{
      author={Borwein, D.},
      author={Maitland, C.},
      author={Skerritt, M.},
       title={Computation of an improved lower bound to {G}iuga's primality
  conjecture},
        date={2013},
     journal={Integers},
      volume={23},
       pages={A67},
}

\bib{GG}{article}{
      author={Giuga, G.},
       title={Su una presumibile propriet\`a caratteristica dei numeri primi},
        date={1950},
     journal={Ist. Lombardo Sci. Lett. Rend. Cl. Sci. Mat. Nat. (3).},
      volume={14},
       pages={511\ndash 528},
}

\bib{GR}{inproceedings}{
      author={Granville, A.},
       title={Arithmetic properties of binomial coefficients i: Binomial
  coefficients modulo prime powers},
        date={2006},
   booktitle={Organic mathematics : proceedings of the organic mathematics
  workshop, december 12-14, 1995. simon fraser university, burnaby, british
  columbia},
      editor={et~al, J.Borwein},
  note={http://www.cecm.sfu.ca/organics/papers/granville/paper/binomial/html/binomial.html},
}

\bib{G}{book}{
      author={Guy, R.~K.},
       title={Unsolved problems in number theory (2nd edition)},
   publisher={Springer},
     address={Berlin},
        date={2004},
}

\bib{K}{misc}{
      author={Kellner, B.C.},
       title={The equivalence of {G}iuga's and {A}goh's conjectures},
        date={2004},
        note={math/0409259},
}

\bib{LPS}{article}{
      author={Luca, F.},
      author={Pomerance, C.},
      author={Shparlinski, I.},
       title={On {G}iuga numbers},
        date={2009},
     journal={Int. J. Mod. Math.},
      volume={4},
       pages={13\ndash 18},
}

\bib{NMSE}{misc}{
      author={Nicolas, A.},
       title={Sums of powers below a prime},
        date={2018},
  note={https://math.stackexchange.com/questions/433678/sums-of-powers-below-a-prime},
}

\bib{S}{misc}{
      author={Scott, Brian~M.},
       title={Congruence for stirling number of first kind $s(n,k)$ when $n$ is
  prime},
        date={2018},
  note={https://math.stackexchange.com/questions/1241649/congruence-for-stirling-number-of-first-kind-sn-k-when-n-is-prime},
}

\bib{ZWS}{unpublished}{
      author={Sun, Z.-W.},
       title={Introduction to {B}ernoulli and {E}uler polynomials},
        date={2002},
        note={Lecture given in Taiwan http://maths.nju.edu.cn/~zwsun/BerE.pdf},
}

\bib{T}{article}{
      author={Tipu, V.},
       title={A note on {G}iuga's conjecture},
        date={2007},
     journal={Canad. Math. Bull.},
      volume={50},
       pages={158\ndash 160},
}

\end{biblist}
\end{bibdiv}
\end{document}